%% file: main.tex
\documentclass[11pt,english,letterpaper]{article}

\usepackage[margin=3cm]{geometry}

\usepackage[english]{babel}
\usepackage[fixlanguage]{babelbib}
\usepackage[utf8]{inputenc}
\usepackage[T1]{fontenc}

\usepackage{acronym} 

\usepackage{lmodern} 

\usepackage[colorlinks]{hyperref} 

\usepackage{xcolor} 
\hypersetup{
  colorlinks=false,
  linkcolor=cyan,
  filecolor=red,      
  urlcolor=blue,
  citecolor=green,
  pdftitle={Orthogonal Projection of Convex Sets with a Differentiable Boundary},
  citebordercolor=green,
  filebordercolor=red,
  linkbordercolor=cyan
}

\usepackage{pgfplots}  
\pgfplotsset{compat = newest}
\pgfdeclarelayer{bg}    
\pgfsetlayers{bg,main}  

\usepackage{amssymb} 
\usepackage{amsmath}
 
\usepackage{derivative} 

\usepackage{tikz}
\usetikzlibrary{arrows.meta} 
\usetikzlibrary{calc,quotes,angles}

\usepackage{amsthm} 
\usepackage{framed}
\usepackage{mdframed}

\newtheorem{Property}{Property}[section]
\newtheorem{Remark}{Remark}[section]
\newtheorem{Theorem}{Theorem}[section]
\newtheorem{Lemma}{Lemma}[section]
\newtheorem{Corollary}{Corollary}[section]

\usepackage{graphicx}
\usepackage{epstopdf}
\usepackage{subcaption} 
\usepackage{stmaryrd}

\usepackage[mathscr]{euscript} 

\usepackage{tcolorbox}
\usepackage{mdframed}

\usepackage{authblk} 

\DeclareMathOperator{\conv}{conv} 
\DeclareMathOperator{\spn}{span} 
\DeclareMathOperator{\Ker}{Ker} 
\DeclareMathOperator{\intr}{int} 
\DeclareMathOperator{\cl}{cl} 
\DeclareMathOperator{\convex}{conv} 

\title{Orthogonal Projection of Convex Sets with a Differentiable Boundary}

\author[1]{Gustave Bainier }
\author[1]{Benoît Marx}
\author[1]{Jean-Christophe Ponsart}

\affil[1]{\small Université de Lorraine, CNRS, CRAN, F-54000 Nancy, France}

\ifpdf
\hypersetup{ pdftitle={Orthogonal Projection of Convex Sets with a Differentiable Boundary} }
\fi

\providecommand{\keywords}[1]
{
  \small	
  \textbf{\textit{Keywords---}} #1
}

\providecommand{\MSCcodes}[1]
{
  \small	
  \textbf{\textit{MSC codes--}} #1
}

\date{\vspace{-5ex}}

\begin{document}
\maketitle

\hrulefill

\begin{abstract}
Given an Euclidean space, this paper elucidates the topological link between the partial derivatives of the Minkowski functional associated to a set (assumed to be compact, convex, with a differentiable boundary and a non-empty interior) and the boundary of its orthogonal projection onto the linear subspaces of the Euclidean space. A system of equations for these orthogonal projections is derived from this topological link. This result is illustrated by the projection of the unit ball of norm $4$ in $\mathbb{R}^3$ on a plane.
\end{abstract}

\\

\MSCcodes{52A20, 53A07}

\hrulefill 

\section{Introduction}

\begin{figure}[h]
\centering
\begin{tikzpicture} 
\def\a{1} \def\b{4}
\def\s{1.5}
\foreach \r in {0,...,8}
{
	\draw[line width=0.5pt,lightgray] plot[samples=100,domain=0:360,smooth,variable=\t] ({\b*(1-(\r / 8))*cos(\t) / \s},{\b*(\r / 8)*sin(\t) / \s});
}
\draw[line width=0.8pt,black] plot[samples=100,domain=0:360,smooth,variable=\t] ({((\b-\a)*cos(\t)+\a*cos((\b-\a)*\t/\a)) / \s},{((\b-\a)*sin(\t)-\a*sin((\b-\a)*\t/\a)) / \s});
\end{tikzpicture}
\caption{The black astroid in the picture above can be seen as the envelope of the family of gray curves $(\mathcal{C}_t)_{t\in(0,1)}$ defined by $\mathcal{C}_t \,:\, \left(x/\left(1-t\right)\right)^{2}+\left(y/t\right)^{2}=1$}
\label{fig:astroid}
\end{figure}
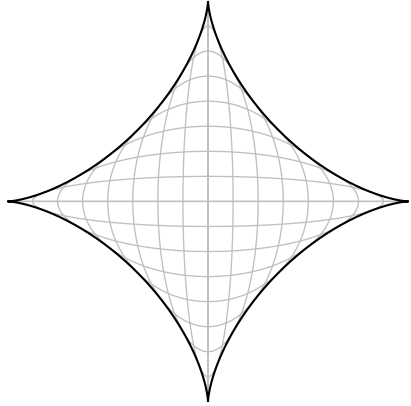

In analytical geometry, given a family of curves $(\mathcal{C}_t)_{t\in \mathbb{R}}$ defined on the plane $\mathbb{R}^2$ by
\begin{equation} 
\mathcal{C}_t : F(x,y,t)=0
\end{equation}
with $F$ a differentiable function, the envelope of $(\mathcal{C}_t)_{t\in \mathbb{R}}$ is defined as the set of points $(x,y) \in \mathbb{R}^2$ such that \cite{eisenhart1909treatise,Pottmannenvelopes}
\begin{equation} 
\exists t\in \mathbb{R},\;\;\;\; F(x,y,t)=0 \;\;\;\;\;\;\frac{\partial F}{\partial t}(x,y,t)=0
\end{equation}
The well-known envelope theorem, mainly used in economics and optimization \cite{Afriat1971,carter2001foundations,Milgrom02envelopeTheorems,Lfgren2011OnET}, provides conditions for the envelope of a family of curves $(\mathcal{C}_t)_{t\in \mathbb{R}}$ to coincide with a single curve tangent to all of the $\mathcal{C}_t$. Under some circumstances, this curve is also the boundary of the region filled by $(\mathcal{C}_t)_{t\in \mathbb{R}}$, and despite this characterization being visually clear (Figure~\ref{fig:astroid}), the authors have not been able to find a satisfying topological discussion on this matter in the literature \cite{Milnor1997-cq,Jottrand2013}.\\

\begin{figure}[h]
\centering
\includegraphics[scale=0.35]{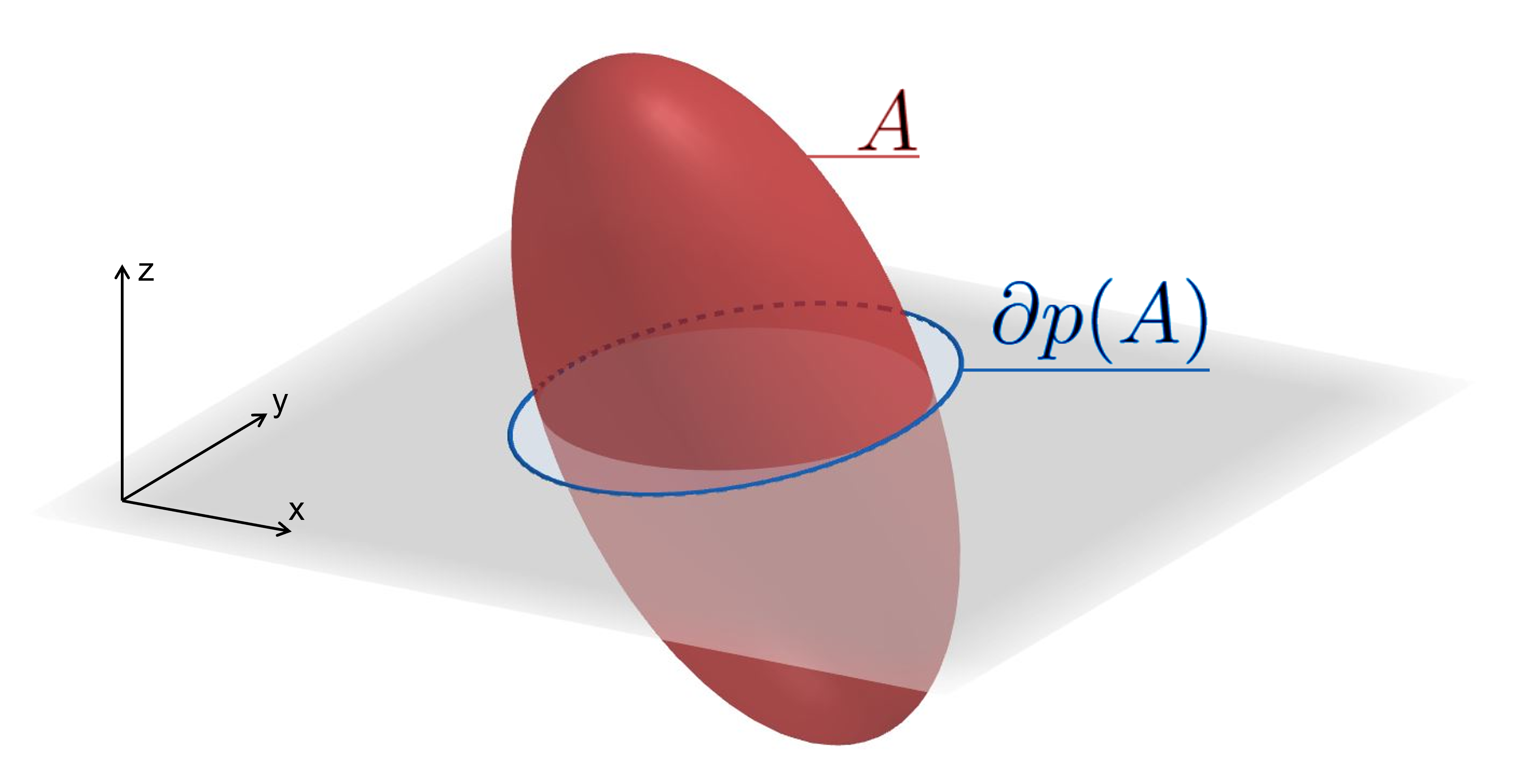}
\caption{$A$, the 3-dimensional ellipsoid in red, is a convex and compact set of $\mathbb{R}^3$. $\partial p(A)$, the boundary of a 2-dimensional ellipsoid with a blue outline, is the boundary of the projection of $A$ along the $z$-axis onto the $xy$-plane represented in gray.}
\label{fig:ellintro}
\end{figure}

Now, given $A$ a convex set of $\mathbb{R}^3$ with a boundary characterized by $F(x,y,z)=0$ where $F$ is differentiable, one can intuitively see by the envelope theorem how characterizing the boundary of $A$ projected along the $z$-axis onto the $xy$-plane relates to the partial derivative of $F$ with respect to $z$ vanishing (Figure~\ref{fig:ellintro}). Moreover, the function $F$ can be obtained from $\mu_A$, the Minkowski functional associated with $A$, usually with the relation $F = \mu_A-1$ \cite{luenberger1968optimization}. In a more general setting, with $E$ a Euclidean space and $A$ a compact and convex set of $E$ with a differentiable boundary and a non-empty interior, the aim of this document is to elucidate the link between the partial derivatives of $\mu_A$ and the boundary of the orthogonal projection of $A$ onto the linear subspaces of $E$. Leveraging results from convex analysis \cite{Rockafellar1970}, a system of equations for the orthogonal projection of $A$ onto any linear subspace of $E$ is obtained. This is the main contribution of the document.\\

The paper is organised as follows: first, in Section~\ref{sec2}, the main definitions and notations used throughout the document are introduced. Then, in Section~\ref{sec3}, preliminary results are derived from topology, convex analysis and properties of the Minkowski functional. These results are applied in Section~\ref{sec4} to elucidate the topological link between the partial derivatives of $\mu_A$ and the boundary of the projection of $A$ onto the linear subspaces of $E$, and a system of equations for the orthogonal projection of $A$ onto the linear subspaces of $E$ is obtained. Section \ref{sec5} provides an illustrative example of the main result of this document by computing the projection of the unit ball of norm $4$ in $\mathbb{R}^3$ on a plane. Finally, Section~\ref{sec6} concludes the document with some application perspectives.

\section{Definitions, Notations}\label{sec2}

$\mathbb{R}$ denotes the field of real numbers. $\mathbb{R}^*$ denotes $\mathbb{R}\setminus \{0\}$. $\mathbb{R}_+$ denotes $[0,+\infty)$. $\mathbb{R}^*_+$ denotes $(0,+\infty)$.\\

Let $E$ denote a Hilbert space of finite dimension over $\mathbb{R}$. $E$ possesses an inner product $\langle\cdot |\cdot\rangle$ which naturally induces a norm $\lVert\cdot\rVert$ and a distance $d(\cdot,\cdot)$ on $E$. $\mathscr{B}_E(x,r)$ denotes the open ball of $E$ centered at $x$ and of radius $r$.\\

Let $A$ and $B$ be two subsets of $E$.  $A+B$ denotes the Minkowski sum of the two sets. $\conv (A)$ denotes the convex hull of $A$ in $E$. $\spn(A)$ denotes the linear span of $A$ in $E$. $\intr_E(A)$, $\cl_E(A)$ and $\partial_E (A)$ denote respectively the interior, the closure and the boundary of $A$ in $E$. $tA$ denotes the scaled set $\{x\in E : x=ty, y\in A\}$. $A$ is said to be absorbing if for all $x\in E$ there exists $t\in\mathbb{R}_+$ such that $x\in tA$.\\

Let $\mathcal{V}$ and $\mathcal{W}$ be two linear subspaces of $E$. $\mathcal{V}\oplus \mathcal{W}$ denotes the direct sum of $\mathcal{V}$ and $\mathcal{W}$. $\mathcal{V}^{\perp}$ denotes the orthogonal complement of $\mathcal{V}$ in $E$. $\dim(\mathcal{V})$ denotes the dimension of $\mathcal{V}$.\\

Let $F$ be another Hilbert space of finite dimension over $\mathbb{R}$ and let $U$ be a subset of $E$. $\mathcal{C}^0(U,F)$ denotes the set of continuous maps from $U$ to $F$. $\mathscr{C}^1(U,F)$ denotes the set of differentiable maps from $U$ to $F$ with a continuous derivative. Given $L$ a linear map from $E$ to $F$, $|||L|||$ denotes the operator norm of $L$. Given $f\in\mathscr{C}^1(U,\mathbb{R})$, $\nabla f (x)$ denotes the gradient of $f$ at $x$.\\ 

The Minkowski functional of $A$ is the map $\mu_A:E \to \mathbb{R}_+$ defined by $\mu_A(x) := \inf \{t \in \mathbb{R}_+^{*} : x \in tA\}$. $A$ has a differentiable boundary if $\mu_A\in \mathscr{C}^1(E\setminus\{0\},\mathbb{R})$.\\

Let $z\in E$, the hyperplane $H$ defined by $H=\Ker(\langle z | \cdot \rangle)$ is called a supporting hyperplane of $A$ at $x\in\partial_E(A)$ if for all $y\in E$, $\mu_A(y)\geq \mu_A(x)+\langle z | y-x \rangle$.\\

In the following, $A$ always denotes a convex, bounded set of $E$ with $0\in \intr_E(A)$ (hence $A$ is absorbing). For all $x\in E$, there exists a unique $x_{\mathcal{V}} \in \mathcal{V}$ and a unique $x_{\mathcal{V}^{\perp}} \in \mathcal{V}^{\perp}$ such that $x = x_{\mathcal{V}}^{\phantom{\perp}} + x_{\mathcal{V}^{\perp}}$. From now on, $p_{\mathcal{V}}$ always denotes the map $x \mapsto x_{\mathcal{V}}$, that is to say the orthogonal projection along $\mathcal{V}^{\perp}$ onto $\mathcal{V}$, with $\mathcal{V}\neq \{0\}$.

\section{Preliminary results}\label{sec3}

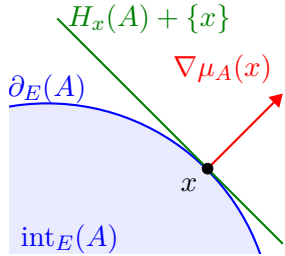
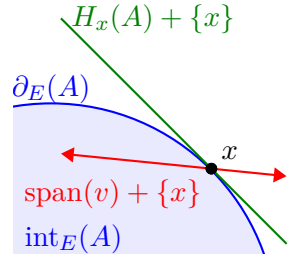
\begin{figure}[h]
\centering
\begin{subfigure}[b]{0.48\textwidth}
\centering
\begin{tikzpicture} 
\begin{scope}
\clip (-1.5,0) rectangle (2.5,3.35);
\draw[line width=0.8pt,blue,fill=blue!40!white,fill opacity=0.2] (-1,-1) circle (3);
\node[text=blue] at (-0.7,0.2) {$\intr_E(A)$};
\node[text=blue] at (-1,2.2) {$\partial_E(A)$};
\draw[line width=0.8pt,black!50!green] (1.1314+1,1.1314-1) -- (1.1314-2,1.1314+2) node[anchor=west,text=black!50!green] {$H_{x}(A)+\{x\}$};
\draw[line width=0.8pt,red,-Triangle] (1.1314,1.1314)--(1.1314+1,1.1314+1) node[anchor=south east,text=red] {$\nabla \mu_A(x)$};
\filldraw[black] (1.1314,1.1314) circle (2pt) node[anchor=north east,text=black]{$x$} ;
\end{scope}
\end{tikzpicture}
\caption{Illustration of Corollary~\ref{firstCorrr}}
\label{subfig01}
\end{subfigure}
\hfill
\begin{subfigure}[b]{0.48\textwidth}
\centering
\begin{tikzpicture} 
\begin{scope}
\clip (-1.5,0) rectangle (2.5,3.35);
\draw[line width=0.8pt,blue,fill=blue!40!white,fill opacity=0.2] (-1,-1) circle (3);
\node[text=blue] at (-0.7,0.2) {$\intr_E(A)$};
\node[text=blue] at (-1,2.2) {$\partial_E(A)$};
\draw[line width=0.8pt,black!50!green] (1.1314+1,1.1314-1) -- (1.1314-2,1.1314+2) node[anchor=west,text=black!50!green] {$H_{x}(A)+\{x\}$};
\draw[line width=0.8pt,red,Triangle-Triangle] (1.1314+1,1.1314-0.1) -- (1.1314-2,1.1314+0.2);
\filldraw[black] (1.1314,1.1314) circle (2pt) node[anchor=south west,text=black]{$x$} ;
\filldraw[black] (1.1314,1.1314) circle (2pt) node[anchor=north east ,text=red] {$\spn(v)+\{x\}$};
\end{scope}
\end{tikzpicture}
\caption{Illustration of Corollary~\ref{iffH}}
\label{subfig02}
\end{subfigure}
\caption{Illustration of the gradient characterization (Corollary~\ref{firstCorrr}) and of the topological characterization (Corollary~\ref{iffH}) of the supporting hyperplane of $A$ at $x \in \partial_E(A)$ when $A$ has a differentiable boundary.}
\label{subfig0}
\end{figure}

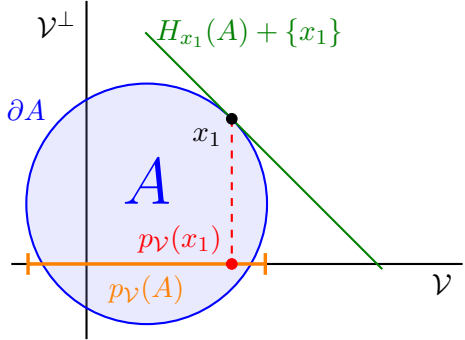
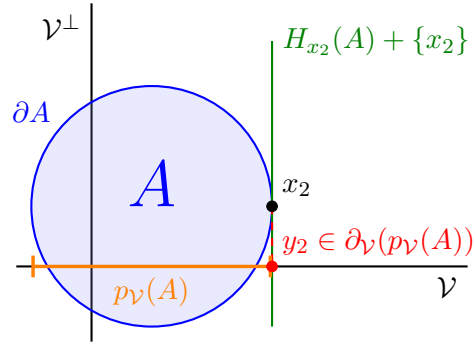
\begin{figure}[h]
\centering
\begin{subfigure}[b]{0.45\textwidth}
\centering
\begin{tikzpicture} 
\draw[line width=0.8pt,black] (-1,0) -- (5,0) node[anchor=north east] {$\mathcal{V}$};
\draw[line width=0.8pt,black] (0,-1) -- (0,3.5) node[anchor=north east] {$\mathcal{V}^{\perp}$};
\draw[line width=0.8pt,blue,fill=blue!40!white,fill opacity=0.2] (1-0.2,1-0.2) circle (1.6);
\node[text=blue] at (0.8,1.1) {\Huge $A$};
\node[text=blue] at (-0.6-0.2,2.05) {$\partial A$};
\draw[line width=0.8pt,black!50!green] (2.1314+2-0.2,2.1314-2-0.2) -- (2.1314-1.15-0.2,2.1314+1.15-0.2) node[anchor=west,text=black!50!green] {$H_{x_1}(A)+\{x_1\}$};
\draw[red,line width=0.8pt,dashed] (2.1314-0.2,2.1314-0.2) -- (2.1314-0.2,0) ;

\draw[|-|,orange,line width=1.3pt] (1-0.2-1.6,0) -- (1-0.2+1.6,0) node [midway, below] {$p_{\mathcal{V}}(A)$};

\filldraw[red] (2.1314-0.2,0) circle (2pt) node[anchor=south east,text=red]{$p_{\mathcal{V}}(x_1)$};
\filldraw[black] (2.1314-0.2,2.1314-0.2) circle (2pt) node[anchor=north east,text=black]{$x_1$} ;
\end{tikzpicture}
\caption{$H_{x_1}(A)$ is the supporting hyperplane of $A$ at $x_1 \in \partial A$. The projection of $x_1$ along $\mathcal{V}^{\perp}$ onto $\mathcal{V}$,  denoted $p_{\mathcal{V}}(x_1)$, is generally unrelated to $H_{x_1}(A)$.}
\label{subfig1}
\end{subfigure}
\hfill
\begin{subfigure}[b]{0.45\textwidth}
\centering
\begin{tikzpicture} 
\draw[line width=0.8pt,black] (-1,0) -- (5,0) node[anchor=north east] {$\mathcal{V}$};
\draw[line width=0.8pt,black] (0,-1) -- (0,3.5) node[anchor=north east] {$\mathcal{V}^{\perp}$};
\draw[line width=0.8pt,blue,fill=blue!40!white,fill opacity=0.2] (1-0.2,1-0.2) circle (1.6);
\node[text=blue] at (0.8,1.1) {\Huge $A$};
\node[text=blue] at (-0.6-0.2,2.05) {$\partial A$};
\draw[line width=0.8pt,black!50!green] (2.4,-0.8) -- (2.4,3) node[anchor=west,text=black!50!green] {$H_{x_2}(A)+\{x_2\}$};
\draw[red,line width=0.8pt,dashed] (2.4,0.8) -- (2.4,0) ;

\draw[|-|,orange,line width=1.3pt] (1-0.2-1.6,0) -- (1-0.2+1.6,0) node [midway, below] {$p_{\mathcal{V}}(A)$};

\filldraw[red] (2.4,0) circle (2pt) node[anchor=south west,text=red]{$y_2 \in \partial_{\mathcal{V}}(p_{\mathcal{V}}(A))$};
\filldraw[black] (2.4,0.8) circle (2pt) node[anchor=south west,text=black]{$x_2$} ;
\end{tikzpicture}
\caption{If $\mathcal{V}^{\perp} \subseteq H_{x_2}(A)$ with $x_2 \in \partial A$, then $y_2 = p_{\mathcal{V}}(x_2)$ is also on the boundary of $p_{\mathcal{V}}(A)$, that is to say $y_2 \in \partial_{\mathcal{V}} (p_{\mathcal{V}}(A))$}
\label{subfig2}
\end{subfigure}
\caption{Illustration of the supporting hyperplanes relation to the orthogonal projection of a convex shape. This relation is formalized in Lemma \ref{tripleLinked}.}
\label{subfig12}
\end{figure}

As stated in the introduction, the partial derivatives of the equation of the boundary of $A$ (a notion of convex analysis) are related to the boundary of the orthogonal projection of $A$ onto the linear subspaces of $E$ (a topological consideration). The main purpose of these preliminary results is to draw a link from convex analysis to topology via the Minkowski functionals associated with $A$. In particular, these preliminary results mainly focus on the link between the gradient of $\mu_A$ and a topological characterization of the supporting hyperplanes of $A$ (Corollary~\ref{firstCorrr}, Corollary~\ref{iffH} and Figure~\ref{subfig0}). The topological characterization of the supporting hyperplanes of $A$ then provides a characterization of the boundary of the projection of $A$ onto $\mathcal{V}$ (Lemma~\ref{tripleLinked} and Figure~\ref{subfig12}), which can finally be linked back to the gradient of $\mu_A$.\\

First, the following classical properties on the Minkowski functional are recalled.

\begin{framed}
\begin{Property}
The Minkowski functional $\mu_A$ satisfies:
\begin{enumerate}
\item For all $x\in E$, $0 \leq \mu_A(x) < +\infty$,
\item For all $x\in E$ and $t\in \mathbb{R}_+$, $\mu_A(tx) = t\mu_A(x)$,
\item For all $x_1,x_2 \in E$, $\mu_A(x_1+x_2) \leq \mu_A(x_1) + \mu_A(x_2)$,
\item $\mu_A \in \mathcal{C}^0(E,\mathbb{R}_+)$
\item $\mu_A^{-1}([0,1)) = \intr_E(A)$, $\mu_A^{-1}([0,1]) = \cl_E(A)$, $\mu_A^{-1}(\{1\}) = \partial_E(A)$
\end{enumerate}
\end{Property}
\end{framed}
\begin{proof}
See Lemma 1 at page 131-132 of \cite{luenberger1968optimization}.
\end{proof}

In particular, item 2 and 3 combined provides the fact that $\mu_A$ is a convex function on $E$. Together with item 5, this establishes a first link between convex analysis and topology. Considering that $A$ has a differentiable boundary (i.e. $\mu_A\in \mathscr{C}^1(E\setminus\{0\},\mathbb{R})$), unicity of the supporting hyperplanes of $A$ is demonstrated using the following result from convex analysis.

\begin{framed}
\begin{Property} \label{subdiff}
Let $f\in \mathscr{C}^{1}(U,\mathbb{R})$ be a convex function. For all $x\in E$, we have:
\begin{equation}\label{subdiffff}
\{z\in E \,:\, \forall y\in U,\, f(y)\geq f(x)+\langle z | y-x \rangle\} = \{\nabla f (x)\}
\end{equation}
\end{Property}
\end{framed}
\begin{proof}
See Theorem 25.1 at page 242 of \cite{Rockafellar1970}.
\end{proof}
\begin{Remark}
The set on the left-hand side of (\ref{subdiffff}) contains the subgradients of $f$ at $x$ and is not necessarily a singleton when $f$ is not differentiable at $x$.
\end{Remark}

Indeed, if $A$ has a differentiable boundary, then $\mu_A$ is a $\mathscr{C}^1$ convex function, and if for all $y\in E$, $\mu_A(y)\geq \mu_A(x)+\langle z | y-x \rangle$, then $\Ker(\langle z | \cdot \rangle)$ is by definition a supporting hyperplane of $A$ at $x\in\partial_E(A)$. Unicity of the supporting hyperplanes of $A$ is obtained from the unicity of such $z$. This links the supporting hyperplanes of $A$ with the gradient of $\mu_A$ (Figure~\ref{subfig01}).

\begin{framed}
\begin{Corollary}[The gradient characterization]\label{firstCorrr}
If $A$ has a differentiable boundary, then there is only one supporting hyperplane of $A$ at $x \in \partial_E(A)$: it is the hyperplane orthogonal to $\nabla \mu_A(x)$. From now on, this supporting hyperplane is denoted $H_x(A)$. Formally, for all $x \in \partial_E(A)$, the following holds:
\begin{equation}
H_x(A)= \Ker(\langle \nabla \mu_A(x) | \cdot \rangle)
\end{equation}
\end{Corollary}
\end{framed}

Now that the supporting hyperplane of $A$ at $x$ is linked with the gradient of $\mu_A$ at $x$, the previous results are now leveraged to obtain a topological characterization of the supporting hyperplanes of $A$. For a convex shape with a differentiable boundary, the supporting hyperplane at a boundary point of this shape is the only hyperplane that, once translated to this point, does not intersect the interior of the shape (Figure~\ref{subfig02}). Lemma~\ref{intersectSame} provides the fact that a supporting hyperplane of $A$ never intersects the interior of $A$, and Lemma~\ref{intersectDiff} provides the fact that, if $A$ has a differentiable boundary, then any affine vector line going through $x \in \partial_E(A)$ that is not included in the supporting hyperplane of $A$ at $x$ will cross the interior of $A$.

\begin{framed}
\begin{Lemma} \label{intersectSame}
If $H$ is a supporting hyperplane of $A$ at $x\in \partial_E(A)$, then $(H+\{x\}) \cap \intr_E(A) = \emptyset$.
\end{Lemma}
\end{framed}
\begin{proof}
By definition of the supporting hyperplane, for all $h\in H$ the following inequality holds $\mu_A(x+h)\geq \mu_A(x)$. Moreover since $x\in\partial_E(A)$, then $\mu_A(x)=1$, which provides $\mu_A(x+h)\geq 1$, hence $( H+\{x\}) \subseteq \mu_A^{-1}([1,+\infty))$, yet $\intr_E(A) = \mu_A^{-1}([0,1))$.
\end{proof}

\begin{framed}
\begin{Lemma}\label{paral}
By parallelism, if $(H+\{x\}) \cap \intr_E(A) = \emptyset$, then $(H+\{x\}) \cap (H+\intr_E(A)) = \emptyset$ as well.
\end{Lemma}
\end{framed}
\begin{proof}
This statement is proved by contraposition.\\
If there exists $y\in (H+\{x\}) \cap (H+\intr_E(A))$, then there exists $z\in\intr_E(A)$ and $h_1,h_2\in H$ such that $y = x+h_1=z+h_2$, providing $z=x+(h_1-h_2)$ where $(h_1-h_2) \in H$, hence $z\in (H+\{x\}) \cap \intr_E(A)$.
\end{proof}

\begin{framed}
\begin{Lemma} \label{intersectDiff}
Suppose $A$ has a differentiable boundary. If $v\notin H_x(A)$, then $(\spn(v)+\{x\}) \cap \intr_E(A) \neq \emptyset$.
\end{Lemma}
\end{framed}
\begin{proof}
This statement is proved by contraposition.\\
Suppose $(\spn(v)+\{x\}) \cap \intr_E(A) = \emptyset$ and consider the function $\phi(t)=\mu_A(x+tv)$. Since $\mu_A\in \mathscr{C}^{1}(E\setminus\{0\},\mathbb{R}_+)$ is a convex function, then $\phi \in \mathscr{C}^{1}(\mathbb{R},\mathbb{R}_+)$ is convex as well. Moreover, since $(\spn(v)+\{x\}) \cap \intr_E(A) = \emptyset$, then for all $t\in \mathbb{R}$, $\phi(t)\geq 1$. Yet $\phi(0)=1$. $t=0$ is therefore a minimum for $\phi$, which implies $\phi'(0)=0$. However, $\phi'(0) = \langle \nabla \mu_A(x) | v \rangle$, hence $v\in\Ker(\langle \nabla \mu_A(x) | \cdot \rangle )$.
\end{proof}

From the Lemmas~\ref{intersectSame} and \ref{intersectDiff}, the following necessary and sufficient condition can be stated, providing a topological characterization of supporting hyperplanes (Figure~\ref{subfig02}) on top of their analytical one (obtained in Corollary~\ref{firstCorrr}):

\begin{framed}
\begin{Corollary}[The topological characterization]\label{iffH}
If $A$ has a differentiable boundary, then $H_x(A)$ contains exactly the directions coming from $x$ that never intersect the interior of $A$. Formally, for all $x \in \partial_E(A)$, the following holds:
\begin{equation}
H_x(A) = \left\{v\in E \,:\, (\spn(v)+\{x\}) \cap \intr_E(A) = \emptyset\right\}
\end{equation}
\end{Corollary}
\end{framed}

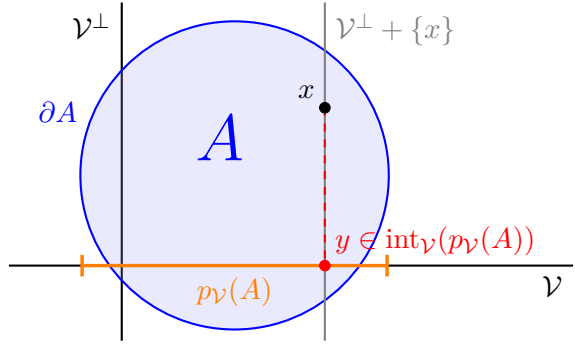
\begin{figure}
\centering
\begin{tikzpicture} 
\draw[line width=0.8pt,black] (-1.5,0) -- (6,0) node[anchor=north east] {$\mathcal{V}$};
\draw[line width=0.8pt,black] (0,-1) -- (0,3.5) node[anchor=north east] {$\mathcal{V}^{\perp}$};
\draw[line width=0.8pt,blue,fill=blue!40!white,fill opacity=0.2] (1.5,1.2) circle (2.05);
\node[text=blue] at (1.3,1.7) {\Huge $A$};
\node[text=blue] at (-0.6-0.25,2.05) {$\partial A$};

\draw[line width=0.8pt,gray] (2.6+0.1,-1) -- (2.6+0.1,3.5) node[anchor=north west] {$\mathcal{V}^{\perp}+\{x\}$};

\draw[red,line width=0.8pt,dashed] (2.6+0.1,2.1) -- (2.6+0.1,0) ;

\draw[|-|,orange,line width=1.3pt] (1.5-2.05,0) -- (1.5+2.05,0) node [midway, below] {$p_{\mathcal{V}}(A)$};

\filldraw[red] (2.6+0.1,0) circle (2pt) node[anchor=south west,text=red]{$y\in \intr_{\mathcal{V}}(p_{\mathcal{V}}(A))$};
\filldraw[black] (2.6+0.1,2.1) circle (2pt) node[anchor=south east,text=black]{$x$} ;
\end{tikzpicture}
\caption{If $x$ is in the interior of $A$, then $y =  p_{\mathcal{V}}(x)$ is in the interior of $p_{\mathcal{V}}(A)$, that is to say $y\in \intr_{\mathcal{V}}(p_{\mathcal{V}}(A))$, and $\mathcal{V}^{\perp}+\{x\}$ crosses the boundary of $A$ multiple times.}
\label{subfig3}
\end{figure}

Before linking the topological characterization of the supporting hyperplanes of $A$ with the boundary of the projection $p_{\mathcal{V}}(A)$ of $A$ onto $\mathcal{V}$, two topological results on the orthogonal projection of $A$ are stated. The first one simply states that the interior of the projection of $A$ is the projection of the interior of $A$ (Lemma~\ref{imgIntr}). The second one states that the projection of the closure of $A$ is also the projection of the boundary of $A$ (Lemma~\ref{imgBndr}). Both are easy to understand visually with the help of Figure~\ref{subfig3}.

\begin{framed}
\begin{Lemma} \label{imgIntr}
If $A$ has a differentiable boundary, then $p_{\mathcal{V}}(\intr_E(A)) = \intr_{\mathcal{V}}(p_{\mathcal{V}}(A))$.
\end{Lemma}
\end{framed}
\begin{proof}This statement is proved by double inclusion.\\

$\subseteq$ This inclusion is a direct consequence of $p_{\mathcal{V}}$ being an open map from $E$ to $\mathcal{V}$. \\

$\supseteq$ Let $y\in\intr_{\mathcal{V}}(p_{\mathcal{V}}(A))$, and $x\in A$ such that $p_{\mathcal{V}}(x)=y$. If $x\in\intr_E(A)$ there is nothing to prove. If $x\in\partial_E(A)$, the following will show by contradiction that $(\mathcal{V}^{\perp}+\{x\})\cap \intr_E(A) \neq \emptyset$, which, thanks to Lemma \ref{intersectDiff}, is equivalent to the existence of $v\in \mathcal{V}^{\perp}$ such that $v\notin H_x(A)$.\\
By contradiction, it is assumed that $\mathcal{V}^{\perp} \subseteq H_x(A)$. By the hyperplane separation theorem, $A$ is contained on one side of $H_x(A)+\{x\}$, hence there is $v\in H_x(A)^{\perp} \setminus \{0\}$ such that for all $t\in \mathbb{R}^*_+$, $x+tv \notin A+ H_x(A)$, therefore $x+tv \notin A+ \mathcal{V}^{\perp}$, and finally $p_{\mathcal{V}}(x+tv) \notin p_{\mathcal{V}}(A)$. However, since $\mathcal{V}^{\perp} \subseteq H_x(A)$, then $v \in \mathcal{V}$, hence $p_{\mathcal{V}}(x+t v)=y+tv$. Since $y\in\intr_{ \mathcal{V}}(p_{ \mathcal{V}}(A))$, by definition of the interior there exists $\delta\in\mathbb{R}^*_+$ such that $\mathscr{B}_{ \mathcal{V}}(y,\delta) \subseteq p_{ \mathcal{V}}(A)$, so in particular there exists $\epsilon \in (0,\delta)$ such that $p_{\mathcal{V}}(x+\epsilon v)= y+\epsilon v \in p_{\mathcal{V}}(A)$, which contradicts that for all $t\in \mathbb{R}^*_+$, $p_{\mathcal{V}}(x+tv) \notin p_{\mathcal{V}}(A)$. Finally $(\mathcal{V}^{\perp}+\{x\})\cap \intr_E(A) \neq \emptyset$.
\end{proof}

\begin{framed}
\begin{Lemma} \label{imgBndr} The following equality holds: $p_{\mathcal{V}}(\cl_E (A)) = p_{\mathcal{V}}(\partial_E (A))$.
\end{Lemma}
\end{framed}
\begin{proof}This statement is proved by double inclusion.\\

 $\subseteq$ Let $y\in p_{\mathcal{V}}(\cl_E (A))$, and $x\in \cl_E (A)$ such that $y=p_{\mathcal{V}}(x)$. If $x\in \partial_{E}(A)$ there is nothing to prove. If $x\in \intr_{E}(A)$, by definition of the interior there exists $\epsilon\in\mathbb{R}^*_+$ such that $\mathscr{B}_E(x,\epsilon) \subseteq A$. Let $v\in \left(\mathscr{B}_E(0,\epsilon) \cap  {\mathcal{V}}^{\perp}\right)\setminus\{0\}$, which guarantees $x+v\in \intr_E(A)$. Since $A$ is bounded, there exists $t\in(1,+\infty)$ such that $x+tv \notin \cl_E(A)$.  Considering the Minkowski functional $\mu_{A+\{-x\}}$,  $x+v\in\intr_E(A)$ translates to $\mu_{A+\{-x\}}(v)<1$, and $x+tv\notin \cl_E(A)$ translates to $\mu_{A+\{-x\}}(tv)>1$. By continuity of  $\mu_{A+\{-x\}}$ the intermediate value theorem provides the existence of $t^*\in(1,t)$ such that $\mu_{A+\{-x\}}(t^*v)=1$, hence $x+t^*v \in \partial_E(A)$. Moreover $x+t^*v \in {\mathcal{V}}^{\perp}+\{y\}$, meaning $p_{\mathcal{V}}(x+t^*v)=y$ (Figure~\ref{subfig3}).\\
 
$\supseteq$ This inclusion is a direct consequence of the inclusion $\partial_E (A) \subseteq \cl_E (A)$.
\end{proof}

With the help of the previous results, the supporting hyperplanes relation to the boundary of the orthogonal projection of $A$ onto $\mathcal{V}$ can be formally stated. Intuitively, when $y\in\mathcal{V}$ is at the boundary of $p_{\mathcal{V}}(A)$, the supporting hyperplane at the pre-image of $y$ by $p_{\mathcal{V}}$ includes ${\mathcal{V}}^{\perp}$, the direction of the projection. Reciprocally, when there is such an alignment, that is to say when ${\mathcal{V}}^{\perp}$ is contained in the supporting hyperplane of the pre-image of $y$ by $p_{\mathcal{V}}$, then $y\in\mathcal{V}$ is at the boundary of $p_{\mathcal{V}}(A)$  (see Figure~\ref{subfig12}). More exactly, the following Lemma holds.

\begin{framed}
\begin{Lemma}\label{tripleLinked}
Let $A$ be closed and have a differentiable boundary. If $y\in p_{\mathcal{V}}(A)$, then the following statements are equivalent:
\begin{enumerate}
  \item $y\in \partial_{\mathcal{V}}(p_{\mathcal{V}}(A))$
  \item $\left\{x\in \partial_E(A) : p_{\mathcal{V}}(x) = y\right\}$ is convex
  \item $\exists x \in \partial_E (A) \,|\,  \begin{cases} p_{\mathcal{V}}(x)=y \\ {\mathcal{V}}^{\perp} \subseteq H_x(A)  \end{cases}$
\end{enumerate}
\end{Lemma}
\end{framed}
\begin{proof} This statement is proved by a circular chain of implications.\\
The notation $B = \left\{x\in \partial_E(A) : p_{\mathcal{V}}(x) = y\right\}$ is used in this proof as a shorthand.\\

$(1) \Rightarrow (2)$ This implication is proved by contraposition.\\
Suppose $B$ is not empty and not convex, hence there exists $z\in \convex(B)\setminus B$. Clearly $\partial_E(A) \subseteq \cl_E(A)$, and the following inclusion is easily checked:
\begin{equation}
    \convex(B) \subseteq \convex\left\{x\in \cl_E(A) : p_{\mathcal{V}}(x) = y\right\}
\end{equation}
moreover, taking $x_1,x_2 \in \left\{x\in \cl_E(A) : p_{\mathcal{V}}(x) = y\right\}$, by linearity of $p_{\mathcal{V}}$, we have for all $t\in[0,1]$, $(tx_1+(1-t)x_2)\in \left\{x\in \cl_E(A) : p_{\mathcal{V}}(x) = y\right\}$, which finally provides the convexity of $\left\{x\in \cl_E(A) : p_{\mathcal{V}}(x) = y\right\}$, hence:
\begin{equation}
    \convex(B) \subseteq \left\{x\in \cl_E(A) : p_{\mathcal{V}}(x) = y\right\}
\end{equation}
This provides the following:
\begin{equation}
\begin{aligned}
\conv(B)\setminus B&\subseteq \left\{x\in \cl_E(A) : p_{\mathcal{V}}(x) = y\right\} \setminus \left\{x\in \partial_E(A) : p_{\mathcal{V}}(x) = y\right\} \\
&\subseteq\left\{x\in \intr_E(A) : p_{\mathcal{V}}(x) = y\right\} 
\end{aligned}
\end{equation}
This provides $z\in \intr_E(A)$ with $p_{\mathcal{V}}(z)=y$. By definition of the interior, there exists $\epsilon \in \mathbb{R}^*_+$ such that $\mathscr{B}_E(z,\epsilon) \subseteq A$. For all $h\in \mathscr{B}_E(0,\epsilon)$, $p_{\mathcal{V}}(z+h) =y+p_{\mathcal{V}}(h)$, and since $|||p_{\mathcal{V}}|||= 1$, then $p_{\mathcal{V}}(h) \in \mathscr{B}_{\mathcal{V}}(0,\epsilon)$, hence $p_{\mathcal{V}}(\mathscr{B}_E(z,\epsilon)) \subseteq \mathscr{B}_{\mathcal{V}}(y,\epsilon) \subseteq p_{\mathcal{V}}(A)$. This finally provides $y\in\intr_{\mathcal{V}}(p_{\mathcal{V}}(A))$.\\ 

$(2) \Rightarrow (3)$ Since $A$ is closed, then, by Lemma \ref{imgBndr}, $y\in p_{\mathcal{V}}(\partial_E(A))$, hence $B\neq \emptyset$. Let $x\in B$ and $v\in {\mathcal{V}}^{\perp}$. The following will show by contradiction that $t\in \mathbb{R}$, $x + t v \notin \intr_E(A)$.\\
Suppose without loss of generality that there exists  $t \in \mathbb{R}^*_+$ such that $x+t v \in \intr_E(A)$. Since $A$ is bounded, with the help of the intermediate value theorem (similarly to Lemma \ref{imgBndr}), there exists $t^*\in(1,+\infty)$ such that $x+ t^* t v \in \partial_E(A)$. This provides $x\in B$, $x + t^* t v \in B$, and $x+ t v \notin B$, yet $B$ should be convex, so there is a contradiction (Figure~\ref{subfig3}). This provides $\spn(v)\cap\intr_E(A) = \emptyset$, hence by Corollary \ref{iffH}, $v\in H_x(A)$.\\

$(3) \Rightarrow (1)$ Let $x\in\partial_E(A)$ be such that $p_{\mathcal{V}}(x)=y$ and ${\mathcal{V}}^{\perp}\subseteq H_x(A)$. Lemma \ref{paral} provides $(H_x(A)+\{x\}) \cap (H_x(A)+\intr_E(A)) = \emptyset$, hence $({\mathcal{V}}^{\perp}+\{x\}) \cap ({\mathcal{V}}^{\perp}+\intr_E(A)) = \emptyset$. Moreover the following equalities hold:
\begin{equation}
\begin{aligned}
({\mathcal{V}}^{\perp}+\{x\}) \cap ({\mathcal{V}}^{\perp}+\intr_E(A))&= p_{\mathcal{V}}^{-1}(\{y\}) \cap p_{\mathcal{V}}^{-1} (p_{\mathcal{V}}(\intr_E(A))) &\\
&= p_{\mathcal{V}}^{-1}(\{y\}) \cap p_{\mathcal{V}}^{-1}(\intr_{\mathcal{V}}(p_{\mathcal{V}}(A))) &[\mbox{Lemma } \ref{imgIntr}] \\
&= p_{\mathcal{V}}^{-1}(\{y\} \cap\intr_{\mathcal{V}}(p_{\mathcal{V}}(A))) & \\
\end{aligned}
\end{equation}
Hence $\{y\} \cap\intr_{\mathcal{V}}(p_{\mathcal{V}}(A)) = \emptyset$, that is to say $y\notin \intr_{\mathcal{V}}(p_{\mathcal{V}}(A))$, providing $y\in \partial_{\mathcal{V}}(p_{\mathcal{V}}(A))$.
\end{proof}

\begin{figure}
\centering
\begin{tikzpicture} 
\node[text=blue] at (2.7,2.2) {$A = \bigcup_{t \in [0,1]} tA$};
\begin{scope}
\clip (-2.05,0) rectangle (2.05,3);
\foreach \i in {0,...,10}
{
	\draw[line width=0.5pt,blue,fill=blue!40!white,fill opacity=0.2] (-\i*0.05,0.7-\i*0.07) circle (2.05-\i*0.2);
}
\end{scope}
\draw[line width=0.8pt,black] (-0.5,0) -- (-0.5,3.5) node[anchor=north east] {$\mathcal{V}^{\perp}$};
\draw[line width=0.8pt,black] (-4,0) -- (4,0) node[anchor=north east] {$\mathcal{V}$};
\clip (-3,0) rectangle (3,-1.5);
\foreach \i in {0,...,10}
{
	
	\draw[line width=0.5pt,blue, opacity=0.2] (-\i*0.05,0.7-\i*0.07) circle (2.05-\i*0.2);
}
\draw[line width=0.8pt,black, opacity=0.2] (-0.5,0) -- (-0.5,-3);
\draw[|-|,orange!95!black,line width=1.5pt] (-2.05,0) -- (2.05,0)  node [midway, below=2pt] {$p_{\mathcal{V}}(A) = \bigcup_{t \in [0,1]}  \partial_{\mathcal{V}} (p_{\mathcal{V}}(tA))$};

\end{tikzpicture}
\caption{Illustration of Lemma~\ref{bigcupthing}, where $A$ is assumed to be closed}
\label{subfig4}
\end{figure}
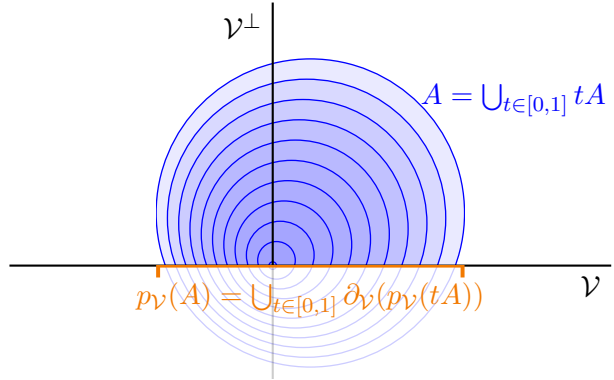

Lastly, the projection of $A$ onto $\mathcal{V}$ can be seen as the union of the boundaries of the projection of $tA$ onto $\mathcal{V}$ with $t\in[0,1]$ (see Figure~\ref{subfig4}). In the next section, the following Lemma will provide a way to go from a statement on the boundary of the projection to a statement on the whole projection $p_{\mathcal{V}}(A)$.

\begin{framed}
\begin{Lemma} \label{bigcupthing} The following equality holds: $p_{\mathcal{V}}(\cl_E(A)) = \bigcup_{t \in [0,1]}  \partial_{\mathcal{V}} (p_{\mathcal{V}}(\cl_E (tA))) $
\end{Lemma}
\end{framed}
\begin{proof}
$\mu_{p_{\mathcal{V}}(\cl_E(A))}$ denotes the Minkowski functional of $p_{\mathcal{V}}(\cl_E(A))$ defined over $\mathcal{V}$. The following equalities hold:
\begin{equation}
\begin{aligned}
p_{\mathcal{V}}(\cl_E(A)) &= \cl_{\mathcal{V}}(p_{\mathcal{V}}(\cl_E(A))) & [p_{\mathcal{V}} \mbox{ continuous} ]\\
&= \mu_{p_{\mathcal{V}}(\cl_E(A))}^{-1}([0,1]) & \\
&= \bigcup_{t \in [0,1]} \mu_{p_{\mathcal{V}}(\cl_E(A))}^{-1}(\{t\}) & \\
&= \bigcup_{t \in [0,1]} \partial_{\mathcal{V}} (t p_{\mathcal{V}}(\cl_E(A))) & \\
&= \bigcup_{t \in [0,1]} \partial_{\mathcal{V}} (p_{\mathcal{V}}(t\cl_E(A))) & [\mbox{linearity of } p_{\mathcal{V}}]\\
p_{\mathcal{V}}(\cl_E(A))  &= \bigcup_{t \in [0,1]} \partial_{\mathcal{V}} (p_{\mathcal{V}}(\cl_E(tA))) & \\
\end{aligned}
\end{equation}
\end{proof}

\section{Characterization of the orthogonal projection of a convex set with a differentiable boundary}\label{sec4}
\input{main_result}



\section{Illustrative example}\label{sec5}

As an illustrative example of Theorem~\ref{mainResult}, this section of the document provides an implicit parametric equation to the projection of the unit ball of norm $4$ of $\mathbb{R}^3$ (denoted $A$) onto the plane $H  : x+y+z = 0$. \\

The Minkowski functional of $A$ is given by

\begin{equation}
\mu_A(x,y,z)=\sqrt[4]{x^4+y^4+z^4}
\end{equation}

After the orthonormal change of basis
\begin{equation}
\left[\begin{array}{c} x \\ y \\ z \end{array}\right] =\left[\begin{array}{ccc} 0 & \sqrt{2/3}  & 1/\sqrt{3} \\  1/\sqrt{2} & -1/\sqrt{6} &  1/\sqrt{3} \\- 1/\sqrt{2} & -1/\sqrt{6} & 1/\sqrt{3} \end{array}\right] \left[\begin{array}{c} u \\ v\\ w\end{array}\right]
\end{equation}
where $w$ is chosen such that $H :w=0$, the function $\eta$ is introduced
\begin{equation}
\eta_A(u,v,w)=\mu_A\left(\sqrt\frac{2}{3}v+\frac{1}{\sqrt{3}}w,\frac{1}{\sqrt{2}}u-\frac{1}{\sqrt{6}}v+\frac{1}{\sqrt{3}}w,-\frac{1}{\sqrt{2}}u-\frac{1}{\sqrt{6}}v+\frac{1}{\sqrt{3}}w\right)
\end{equation}
For all $(u,v,w)\neq (0,0,0)$, its partial derivative with respect to $w$ is given by
\begin{equation}
\begin{aligned}
\frac{\partial \eta_A}{\partial w}(u,v,w)&=\frac{\partial }{\partial w}\left[\sqrt[4]{\left[\sqrt\frac{2}{3}v+\frac{1}{\sqrt{3}}w\right]^4+\left[\frac{1}{\sqrt{2}}u-\frac{1}{\sqrt{6}}v+\frac{1}{\sqrt{3}}w\right]^4+\left[-\frac{1}{\sqrt{2}}u-\frac{1}{\sqrt{6}}v+\frac{1}{\sqrt{3}}w\right]^4}\right] \\
&=\left[\frac{1}{3}w^3+(u^2+v^2)w  -\frac{\sqrt{2}}{2}u^2 v+\frac{\sqrt{2}}{6}v^3 \right]\eta_A^{-3}(u,v,w) \\
\end{aligned}
\end{equation}
Since $\eta_A^{-3}(u,v,w)>0$, studying $w$ such that $\frac{\partial \eta_A}{\partial w}(u,v,w) = 0$ is equivalent to the study of the solutions to the depressed cubic equation 
\begin{equation} \label{polyn3}
X^3+3(u^2+v^2)X+ \frac{\sqrt{2}}{2}v\left(v^2-3u^2\right) = 0
\end{equation}
which discriminant is given by
\begin{equation}
\Delta = -\left(108(u^2+v^2)^3+\frac{27}{2}v^2\left(v^2-3u^2\right)^2\right)
\end{equation}
It is easily verified that $\Delta \leq 0$, hence there is only one real root $w^*$ satisfying \eqref{polyn3}, and it is given by Cardano's formula \cite{Van_der_Waerden2003-dm}
\begin{equation}
w^*=\sqrt[3]{-\frac{\sqrt{2}}{4}v\left(v^2-3u^2\right)-\sqrt{\delta(u,v)}}+\sqrt[3]{-\frac{\sqrt{2}}{4}v\left(v^2-3u^2\right)+\sqrt{\delta(u,v)}}
\end{equation}
where
\begin{equation}
\delta(u,v)=\frac{1}{8}v^2\left(v^2-3u^2\right)^2+(u^2+v^2)^3
\end{equation}
Finally, Theorem~\ref{mainResult} provides that the projection of $A$ onto $H$ is given by the $(u,v)\in\mathbb{R}^2$ satisfying
\begin{equation}
\eta_A\left(u,v,\sqrt[3]{-\frac{\sqrt{2}}{4}v\left(v^2-3u^2\right)-\sqrt{\delta(u,v)}}+\sqrt[3]{-\frac{\sqrt{2}}{4}v\left(v^2-3u^2\right)+\sqrt{\delta(u,v)}}\right) \leq 1
\end{equation}
which is plotted in Figure \ref{fig:example} below.

\begin{figure}[h]
\centering
\includegraphics[scale=0.35]{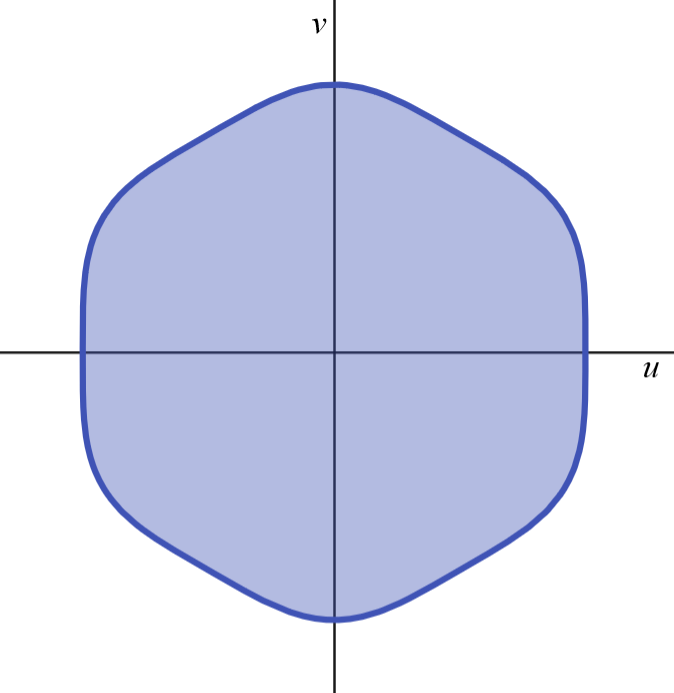}
\caption{Shape of the projection of the unit ball of norm $4$ of $\mathbb{R}^3$ onto the plane $H  : x+y+z = 0$}
\label{fig:example}
\end{figure}

\section{Conclusion}\label{sec6}

In this study, the topological link between the partial derivatives of the Minkowski functional associated with $A$ (a compact and convex set of a Euclidean space $E$) and the boundary of the projection of $A$ onto the linear subspaces of $E$ was elucidated. This topological link provided a system of equations for the orthogonal projection of $A$ onto the linear subspaces of $E$. \\

Some applications of these results can be found for engineering, in particular in fault detection schemes for the diagnosis of dynamical systems. Indeed, model-based fault detection consists in identifying when a fault occurs in a dynamical system by analysing the discrepancies between the system inputs and outputs and their expected values provided by the model \cite{Ding2013}. These discrepancies are generally used to generate residual signals for system diagnosis. However, these signals being subject to the system perturbations and to measurement noises, one of the challenge of fault detection is to distinguish the unavoidable noise from an actual fault in the process \cite{Basseville1993,whalen2013detection}. For example, residuals obtained using a parity space approach to fault detection are generally simply affected by a projection of this noise, hence, knowing a noise bounding shape, an exact threshold to detect a fault could be the boundary of the projection of this bounding shape. 

\addcontentsline{toc}{section}{References}

\bibliographystyle{unsrt}
\bibliography{references}

\end{document}

%% file: main_result.tex
The main result of this document consists in obtaining a system of equations that characterizes the orthogonal projection of the closure of $A$ on a linear subspace $\mathcal{V}\neq \{0\}$ when $A$ has a differentiable boundary. To obtain this system of equations, the following Minkowski functional of two variables is introduced:
\begin{equation}
\begin{aligned}
\eta_A : {\mathcal{V}} \times {\mathcal{V}}^{\perp} &\to \mathbb{R} \\
\left(x_{\mathcal{V}}^{\phantom{\perp}},x_{{\mathcal{V}}^{\perp}}\right) &\mapsto \mu_A\left(x_{\mathcal{V}}^{\phantom{\perp}}+x_{{\mathcal{V}}^{\perp}}\right)
\end{aligned}
\end{equation}
From now on, $\frac{\partial \eta_A}{\partial x_{{\mathcal{V}}}}$ denotes the partial derivative of $\eta_A$ with respect to $x_{{\mathcal{V}}}$ and $\frac{\partial \eta_A}{\partial x_{{\mathcal{V}}^{\perp}}}$ denotes the partial derivative of $\eta_A$ with respect to $x_{{\mathcal{V}}^{\perp}}$.\\

The link between the partial derivatives of $\eta_A$ and the boundary of the orthogonal projection of $A$ onto the linear subspaces of $E$ is explicitly written and leveraged in the proof of this characterization.

\begin{framed}
\begin{Theorem}\label{mainResult}
If $A$ is a compact and convex set of $E$ with a differentiable boundary and $0\in\intr_E(A)$, then, for all projection $p_{\mathcal{V}}$ such that $\mathcal{V}\neq \{0\}$, the following equality holds:
\begin{equation}
p_{\mathcal{V}}(A) = \left\{ x_{\mathcal{V}}^{\phantom{\perp}} \in {\mathcal{V}} \, : \, \exists x_{{\mathcal{V}}^{\perp}}\in {\mathcal{V}}^{\perp} \,|\, \begin{cases} \eta_A\left(x_{\mathcal{V}}^{\phantom{\perp}},x_{{\mathcal{V}}^{\perp}}\right) \leq 1  \\ x_{\mathcal{V}}^{\phantom{\perp}}+x_{{\mathcal{V}}^{\perp}}\neq 0 \Rightarrow \frac{\partial \eta_A}{\partial x_{{\mathcal{V}}^{\perp}}}\left(x_{\mathcal{V}}^{\phantom{\perp}},x_{{\mathcal{V}}^{\perp}}\right) = 0\end{cases} \right\}
\end{equation}
\end{Theorem}
\end{framed}

\begin{proof}
If $t=0$ then $tA=\{0\}=p_{\mathcal{V}}(tA)$, hence for all $x_{{\mathcal{V}}^{\perp}} \in {\mathcal{V}}^{\perp}$, the equality $\eta_A(0,x_{{\mathcal{V}}^{\perp}}) = \mu_{p_{\mathcal{V}}(A)}(0)$ holds, and there is nothing to prove.\\
If $t\in\mathbb{R}^*_+$, thanks to Lemma \ref{tripleLinked}, the following equivalence holds:
\begin{equation}
y\in \partial_{\mathcal{V}}(p_{\mathcal{V}}(tA)) \Leftrightarrow \exists x \in \partial_E (tA) \,|\, \begin{cases} p_{\mathcal{V}}(x)=y \\
{\mathcal{V}}^{\perp} \subseteq H_x(tA) \end{cases}
\end{equation}
For all $x\in \partial_E (tA)$, $H_x(tA)= \Ker(\langle \nabla \mu_{tA}(x) | \cdot \rangle)$, and since $t\neq 0$, then $\nabla \mu_{tA}(x) = \nabla \mu_{A}(x)$. Moreover for all $h \in E$, $x_{\mathcal{V}},h_{\mathcal{V}} \in {\mathcal{V}}$ and $x_{{\mathcal{V}}^{\perp}}, h_{{\mathcal{V}}^{\perp}} \in {\mathcal{V}}^{\perp}$ such that $x = x_{\mathcal{V}}^{\phantom{\perp}}+x_{{\mathcal{V}}^{\perp}}$ and $h = h_{\mathcal{V}}^{\phantom{\perp}}+h_{{\mathcal{V}}^{\perp}}$, the following equality holds:
\begin{equation}
\begin{aligned}
\langle \nabla \mu_{A}(x) | h \rangle 
&= \frac{\partial \eta_A}{\partial x_{{\mathcal{V}}}}\left(x_{\mathcal{V}}^{\phantom{\perp}},x_{{\mathcal{V}}^{\perp}}\right)  h_{\mathcal{V}}^{\phantom{\perp}}   +   \frac{\partial \eta_A}{\partial x_{{\mathcal{V}}^{\perp}}}\left(x_{\mathcal{V}}^{\phantom{\perp}},x_{{\mathcal{V}}^{\perp}}\right)   h_{{\mathcal{V}}^{\perp}}
\end{aligned}
\end{equation}
hence the following equivalences hold:
\begin{equation}
\begin{aligned}
&y\in \partial_{\mathcal{V}}(p_{\mathcal{V}}(tA))\Leftrightarrow \exists  x_{{\mathcal{V}}^{\perp}} \in {\mathcal{V}}^{\perp}\,|\,  \begin{cases}y+x_{{\mathcal{V}}^{\perp}} \in \partial_E (tA) \\ \frac{\partial \eta_A}{\partial x_{{\mathcal{V}}^{\perp}}}(y,x_{{\mathcal{V}}^{\perp}}) = 0 \end{cases} \\
\mbox{i.e. }&y\in \partial_{\mathcal{V}}(p_{\mathcal{V}}(tA)) \Leftrightarrow  \exists  x_{{\mathcal{V}}^{\perp}} \in {\mathcal{V}}^{\perp}\,|\,  \begin{cases}\eta_A(y,x_{{\mathcal{V}}^{\perp}}) = t  \\ \frac{\partial \eta_A}{\partial x_{{\mathcal{V}}^{\perp}}}(y,x_{{\mathcal{V}}^{\perp}}) = 0 \end{cases}
\end{aligned}
\end{equation}
For $t=1$, this last equivalence provides the link between the partial derivatives of $\eta_A$ and the boundary of the orthogonal projection of $A$ onto the linear subspaces of $E$.\\

Finally, Lemma \ref{bigcupthing} provides:
\begin{equation}
\begin{aligned}
p_{\mathcal{V}}(A) &= 
p_{\mathcal{V}}(\cl_E(A)) \\
&= \bigcup_{t \in [0,1]}  \partial_{\mathcal{V}} (p_{\mathcal{V}}(tA)) \\
&=   \{0\}\cup \bigcup_{t \in (0,1]} \left\{ x_{\mathcal{V}}^{\phantom{\perp}} \in {\mathcal{V}} \, : \, \exists x_{{\mathcal{V}}^{\perp}}\in {\mathcal{V}}^{\perp} \,|\, \begin{cases} \eta_A\left(x_{\mathcal{V}}^{\phantom{\perp}},x_{{\mathcal{V}}^{\perp}}\right) = t  \\ \frac{\partial \eta_A}{\partial x_{{\mathcal{V}}^{\perp}}}\left(x_{\mathcal{V}}^{\phantom{\perp}},x_{{\mathcal{V}}^{\perp}}\right) = 0\end{cases} \right\}\\
p_{\mathcal{V}}(A) &= \left\{ x_{\mathcal{V}}^{\phantom{\perp}} \in {\mathcal{V}} \, : \, \exists x_{{\mathcal{V}}^{\perp}}\in {\mathcal{V}}^{\perp} \,|\, \begin{cases} \eta_A\left(x_{\mathcal{V}}^{\phantom{\perp}},x_{{\mathcal{V}}^{\perp}}\right) \leq 1  \\ x_{\mathcal{V}}^{\phantom{\perp}}+x_{{\mathcal{V}}^{\perp}}\neq 0 \Rightarrow \frac{\partial \eta_A}{\partial x_{{\mathcal{V}}^{\perp}}}\left(x_{\mathcal{V}}^{\phantom{\perp}},x_{{\mathcal{V}}^{\perp}}\right) = 0\end{cases} \right\}
\end{aligned}
\end{equation}
which concludes the proof.
\end{proof}

Given a compact and convex set of $E$ with a differentiable boundary and a non-empty interior, there exists a translation so that the origin of $E$ is in the interior of the translated set, hence this new set is absorbing. Given a good translation of $A$, the main result of this document can therefore be extended without difficulty to a more general setting where $A$ simply denotes a compact and convex set of $E$ with a differentiable boundary and a non-empty interior.

\begin{framed}
\begin{Corollary}
Keeping the assumptions of Theorem~\ref{mainResult}, the following equality holds:
\begin{equation} \label{eq:firstCor}
\mu_{p_{\mathcal{V}}(A)}(x) = \begin{cases}\inf \left\{t\in \mathbb{R}^*_+ \, : \, \exists x_{{\mathcal{V}}^{\perp}}\in {\mathcal{V}}^{\perp} \,|\, \begin{cases} \eta_A\left(x,x_{{\mathcal{V}}^{\perp}}\right) \leq t  \\ x+x_{{\mathcal{V}}^{\perp}}\neq 0 \Rightarrow \frac{\partial \eta_A}{\partial x_{{\mathcal{V}}^{\perp}}}\left(x,x_{{\mathcal{V}}^{\perp}}\right) = 0\end{cases} \right\}  &\mbox{ if } x\in \mathcal{V} \\ +\infty &\mbox{ if } x\notin \mathcal{V}\end{cases}
\end{equation}
Moreover, if $V$ and $V^{\perp}$ denote the matrices whose columns are resp. formed by $(v_1,\dots ,v_m)$ a basis to $\mathcal{V}$ and $(v_{m+1}, \dots, v_n)$ a basis to $\mathcal{V}^{\perp}$, then the following equality holds:
\begin{equation}  \label{eq:secondCor}
\mu_{PA}(y) = \mu_{p_{\mathcal{V}}(A)}(Vy)
\end{equation}
with $P=\left[\begin{array}{cc} I_m & 0\end{array}\right]\left[\begin{array}{cc} V & V^{\perp}\end{array}\right]^{-1}$ and where $y\in\mathbb{R}^m$ is expressed in the $(v_1,\dots ,v_m)$ basis.
\end{Corollary}
\end{framed}
\begin{proof}
Equation \eqref{eq:firstCor} is easily derived by replacing the interval $[0,1]$ by the interval $[0,t]$ in the proof of Theorem~\ref{mainResult}. Equation~\eqref{eq:secondCor} is a trivial consequence of $p_{\mathcal{V}}(A) = VPA$.
\end{proof}